\documentclass[a4paper]{amsart}
\usepackage[new]{old-arrows}
\usepackage[T1]{fontenc}
\usepackage[english]{babel}
\usepackage{kpfonts}

\usepackage{amsthm,amsmath,amssymb,amsfonts,euscript,graphicx,hyperref,
indentfirst}

\usepackage{enumerate}
\usepackage[all]{xy}
\usepackage{tikz-cd}
\usepackage{titlepic}
\usetikzlibrary{decorations.pathreplacing,calligraphy}

\usepackage{hyperref}
\usepackage{cleveref}

\usepackage{amssymb,amsfonts}
\usepackage{tikz}

\usetikzlibrary{arrows,snakes,backgrounds}
\usepackage{color}   
\usepackage{marginnote}
\usepackage{tikz-cd}
\usetikzlibrary{positioning}
\usepackage{enumitem}
\usepackage{hyperref}
\usetikzlibrary{arrows}

\usepackage{mathtools}

\DeclarePairedDelimiter\floor{\lfloor}{\rfloor}
\newtheorem{thm}{Theorem}[section]
\newtheorem{cor}[thm]{Corollary}

\newtheorem{prop*}[thm]{Proposition}

\newtheorem{conj}[thm]{Conjecture}
\newtheorem{quest}[thm]{Question}

\newtheorem{problem}[thm]{Problem}

\theoremstyle{definition}
\newtheorem{defn}[thm]{Definition}

\theoremstyle{remark}
\newtheorem{rem}[thm]{Remark}

\newcommand{\bC}{\mathbb{C}}

\newcommand{\bQ}{\mathbb{Q}}
\newcommand{\bR}{\mathbb{R}}
\newcommand{\bS}{\mathbb{S}}

\newcommand{\bZ}{\mathbb{Z}}

\newcommand\Diff{\mathrm{Diff}}
\newcommand\Homeo{\mathrm{Homeo}}

\newcommand\BDiff{\mathrm{BDiff}}

\newcommand{\BcdDiff}[1]{\mathrm{BDiff}_c(#1)^{\delta}}
\newcommand{\cdDiff}[1]{\mathrm{Diff}_c(#1)^{\delta}}

\newcommand{\BdrDiff}[1]{\mathrm{BDiff}^{r}(#1)^{\delta}}
\newcommand{\BcdrDiff}[1]{\mathrm{BDiff}_c^{r}(#1)^{\delta}}

\newcommand{\drDiff}[1]{\mathrm{Diff}^{r}(#1)^{\delta}}
\newcommand{\BdDiffo}[1]{\mathrm{BDiff}_{0}(#1)^{\delta}}
\newcommand{\BdrDiffo}[1]{\mathrm{BDiff}^r_{0}(#1)^{\delta}}
\newcommand{\drDiffo}[1]{\mathrm{Diff}^r_{0}(#1)^{\delta}}

\newcommand{\hcoker}{/\!\!/}

\newcommand{\tH}{\text{\textnormal{Homeo}}}

\newcommand{\BH}{\mathrm{B}\text{\textnormal{Homeo}}}

\newcommand{\tdH}{\text{\textnormal{Homeo}}^{\delta}}

\newcommand{\BdH}{\mathrm{B}\text{\textnormal{Homeo}}^{\delta}}

\titlepic{\includegraphics[width=\textwidth]{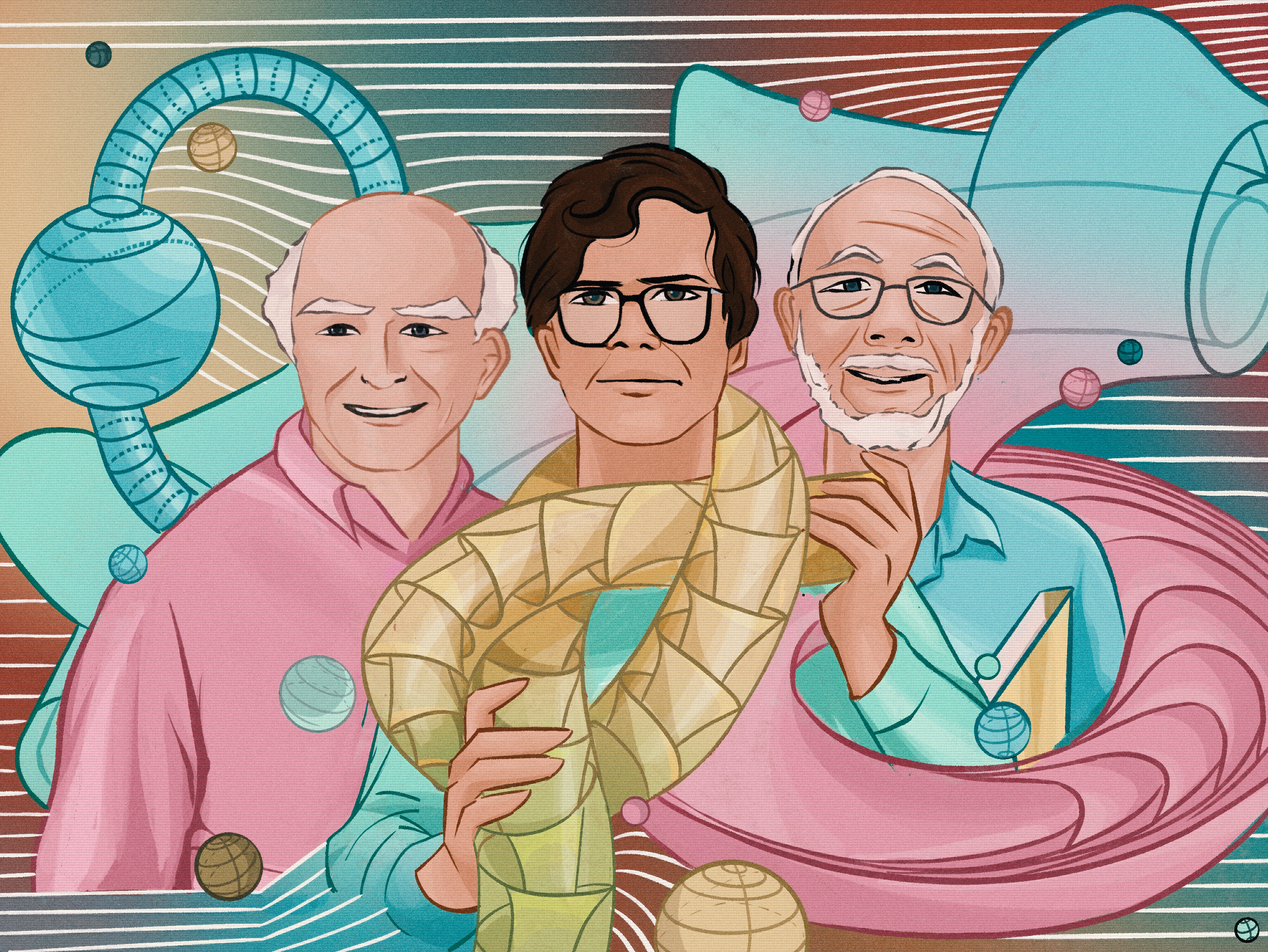}}
\title[Foliations and diffeomorphism groups]{
Foliations and diffeomorphism groups \\ \, \\ {\Small Sam Nariman}
}

\author{Sam Nariman}

\begin{document}

\maketitle
\section{Introduction}
Thurston at his ICM address in 1974 asked ``Given a large supply of some sort of fabrics, what kind of manifolds can be made from it, in a way that the patterns match up along the seams?'' And he alluded that ``For open manifolds, Gromov's theorem gives a good answer for a wide variety of fabrics. The techniques needed to analyze such a question on a closed manifold are usually different, at least to a casual eye''. In this article, we shall focus on one aspect of this difference between closed and open manifolds that is inspired by a property of diffeomorphism groups called {\it fragmentation property}.

But let us first define the basic objects that we would like to study. Recall that manifolds are modeled locally on Euclidean spaces $\bR^n$ along with a subgroup of the homeomorphism group $\tH(\bR^n)$ as the permissible transition functions for a given structure on the manifold. Now a  $C^r$-\emph{foliation} of codimension $k$ on an $n$-dimensional manifold $M$ is determined by a {\it foliated atlas} $\{(U_i, \phi_i)\}_i$ where $\{U_i\}_i$ is an open cover of $M$ and $\phi_i\colon U_i\to \bR^n$ are maps such that for overlapping pairs $U_i$ and $U_j$, the transition function $\phi_j\circ \phi_i^{-1}$ is a $C^r$ diffeomorphism from $\phi_i(U_i\cap U_j)$ 
to $\phi_j(U_i\cap U_j)$ such that the function $\phi_j\circ \phi_i^{-1}$ is of the form
\[
\phi_j\circ \phi_i^{-1}(x,y)=(f(x,y), g(y)),
\]
where $x$ denotes the first $n-k$ coordinates and $y$ denotes the last $k$ coordinates. If $y$ is fixed, then so is $g(y)$, hence the property of having the last $k$ coordinates being fixed is preserved by the change of coordinates. Therefore the sheets that are given by fixing the last $k$ coordinates glue together globally to form leaves of the foliation. In particular, this structure gives a decomposition of $M$ into $(n-k)$-submanifolds (the leaves).

\begin{figure}[h]
  \centering
    \includegraphics[width=0.4\textwidth]{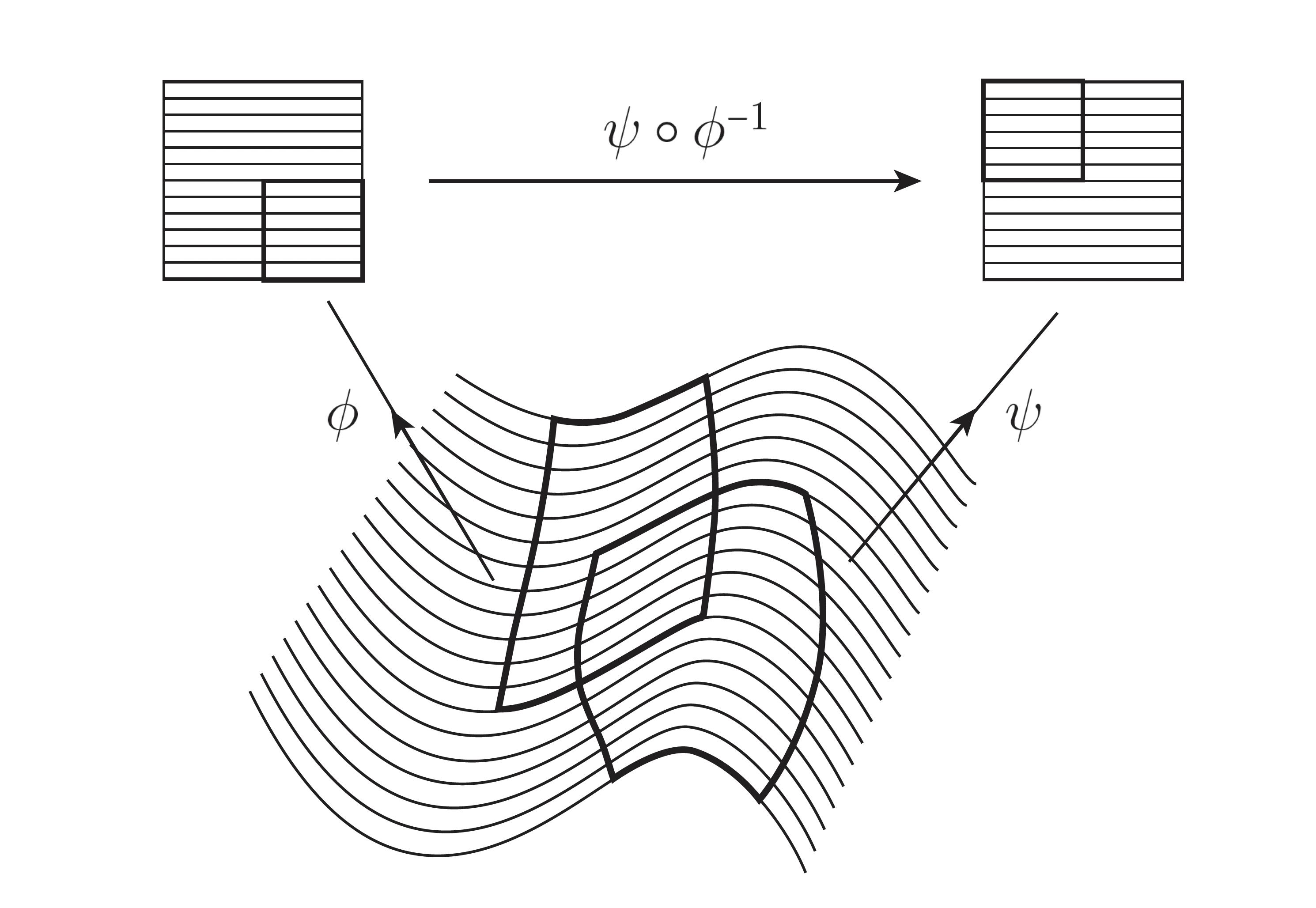}
    \caption{Foliation charts}

\end{figure}
Let us give two examples that will show up later in this article. Consider the foliation of $\bR^2$ by straight parallel lines. Since such foliations are invariant under translations in the $x$-axis and $y$-axis, they induce smooth codimension-one foliations on the $2$-torus $T^2=\bR^2/\bZ^2$ (see Figure 2) and they are called linear foliations on torus. Note that if the angle of the lines with the $x$-axis is irrational, the induced linear foliation has no compact leaf.
\begin{figure}[h]
  \centering
    \includegraphics[width=0.6\textwidth]{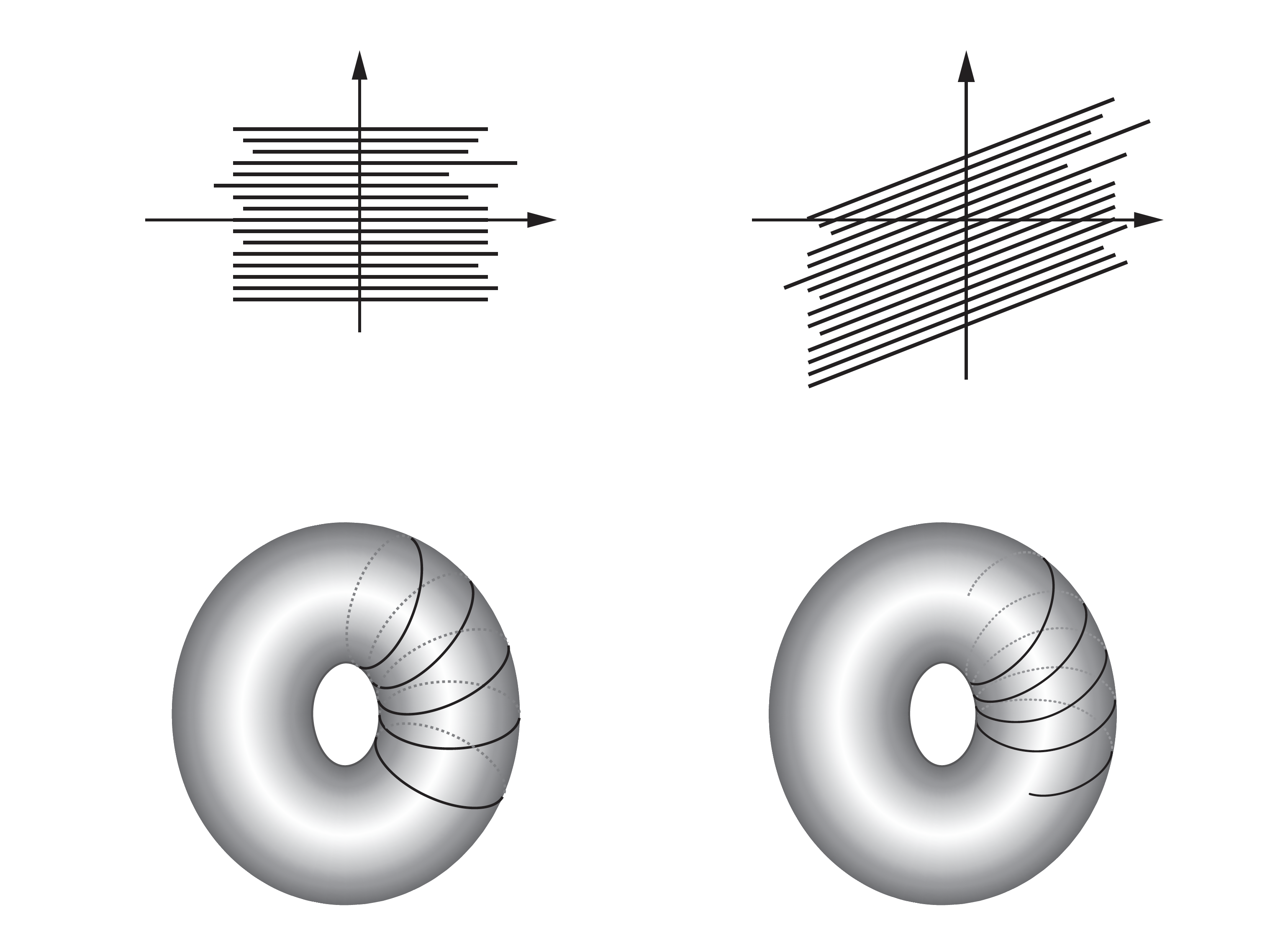}
    \caption{Linear foliations on $2$-torus}

\end{figure}
Another classical example is the Reeb foliation which is a smooth codimension-one foliation on the solid torus. The boundary torus is the only compact leaf and the other leaves are diffeomorphic to $ 2$-planes as in the figure below.
\begin{figure}[h]
  \centering
    \includegraphics[width=0.6\textwidth]{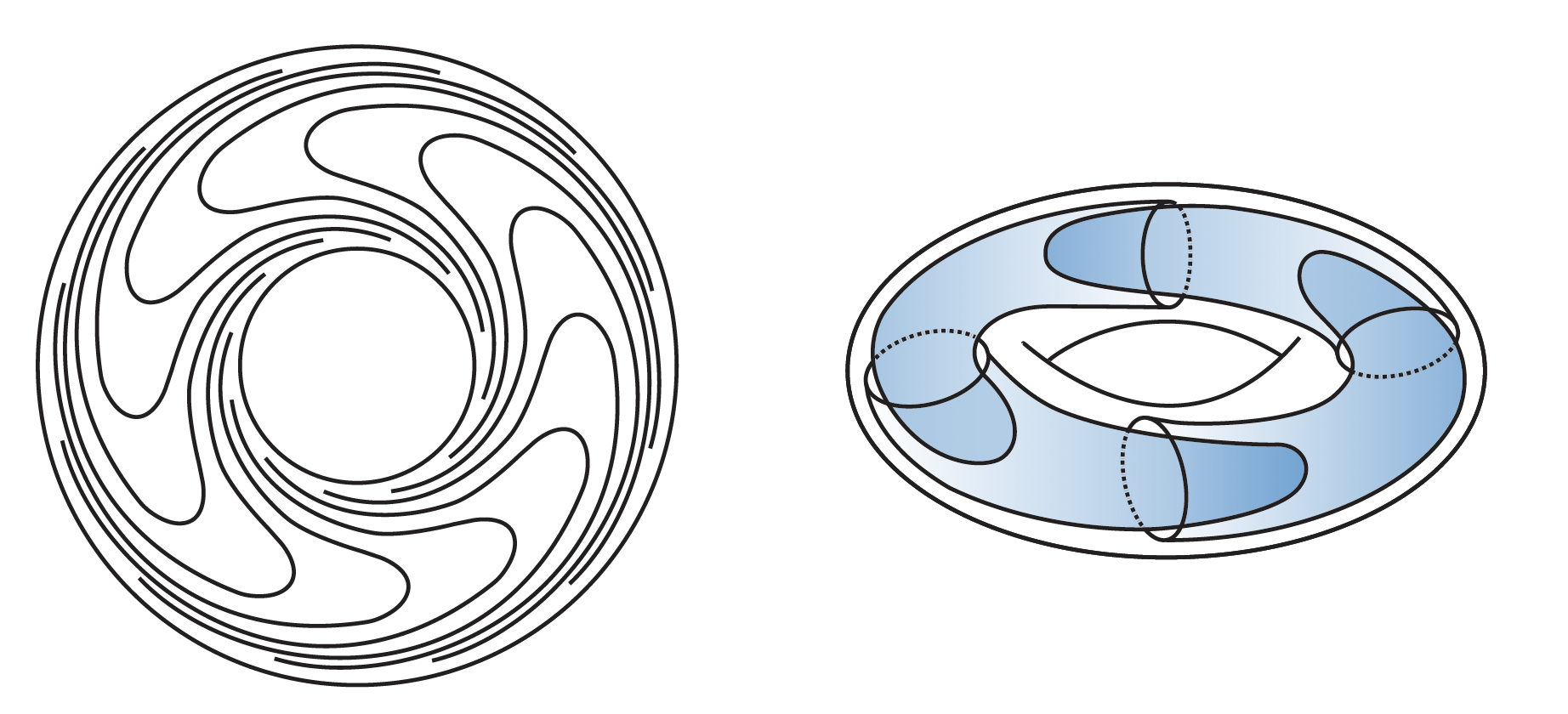}
    \caption{Schematic of $2$-dimensional and $3$-dimensional Reeb foliations on an annulus and on a solid torus}

\end{figure}
One can glue two of these foliated solid tori along the torus boundary to obtain smooth codimension-one foliation on $S^3$.

To admit a codimension $k$ foliation, the manifold $M$ must have a plane field of dimension $n-k$ (the tangents to the leaves). A plane field is called {\it integrable} if it is given by the tangent-plane field to a foliation.   A classical theorem that is usually ascribed to Frobenius but in fact goes back to the earlier work of Clebsch and Deahna (\cite[Theorem 1]{MR0343289}) gives a necessary and sufficient condition for a plane field to be integrable. Sections of a plane field $\xi\subset TM$ are vector fields. The integrability condition for a plane field $\xi$ is that the Lie bracket of $C^1$ sections of $\xi$ should  also be a section of $\xi$. 

We can think of foliations as plane fields satisfying the integrability condition.  To find a plane field of dimension $n-k$ on $M$ is essentially a homotopy theory question. But the integrability condition is of a different type, in local coordinates it can be described as a system of partial differential equations (hence, it is not an open condition). So after finding a plane field, it is more natural to ask whether we can deform this plane field to become integrable in which case we say it is integrable up to homotopy.  Bott found a topological obstruction for plane fields $\xi\subset TM$ to be integrable up to homotopy, which is now known as the Bott vanishing theorem. He showed that  Pontryagin classes of the quotient bundle $TM/\xi$ must vanish in degrees larger than $2\text{dim}(TM/\xi)$. Using this obstruction, for $n>1$ and odd, he showed that each complex codimension-one holomorphic subbundle of $T\bC P^n$ is not integrable up to homotopy. 

``Gromov's theorem'' in Thurston's quote refers to the powerful theory of h-principle (homotopy principle) that Gromov developed to reduce the problem of the existence of solutions of differential relations (differential equations and inequalities) in differential geometry to essentially homotopy theory questions. Philosophically, h-principle theorems relate the topology of a geometric object which is more rigid to a homotopy theoretical object that is more flexible to work with. The first famous example is Smale's sphere eversion which says intuitively that one can turn the sphere $S^2$ in $\mathbb{R}^3$ inside out smoothly without cutting or tearing it or creating any crease. More precisely, the space of immersions of $S^2$ into $\bR^3$ is a connected space. Smale showed that this space of immersions is homotopy equivalent to the space of bundle monomorphisms between $TS^2$ and $T\bR^3$, which can be seen to be connected using algebraic topology. 

Bott's obstruction shows that the existence of foliations cannot be solely reduced to the homotopy theory question of the existence of plane fields.  So we need more homotopical information to reduce the problem of the existence of foliations to a homotopy theory question. Haefliger introduced in \cite{MR0100269} a notion that is now called Haefliger structures on manifolds that are more flexible than foliations in  the following sense and sometimes are called ``singular'' foliations. In order to pull back a foliation $\mathcal{F}$ on $M$ under the map $f\colon N\to M$, the map $f$ has to be transverse to $\mathcal{F}$ (i.e. $df(TN)+T\mathcal{F}= TM$). However, Haefliger structures are more flexible than foliations since they can be pulled back under all continuous maps. Because of this flexibility, one can study them as objects in algebraic topology. 

\begin{defn}A $C^r$ Haefliger structure  $\mathcal{H}$ of codimension $k$ on $M$ is given by the triple $(\nu \mathcal{H}, U, \mathcal{F})$ where $\nu\mathcal{H}$ is a $k$-dimensional $C^r$ vector bundle over $M$
which is called the normal bundle of the Haefliger structure $\mathcal{H}$, the open subset $U\subset \nu\mathcal{H}$ around the zero section, and a $C^r$ foliation of codimension $k$ in $U$ that is transverse to the fibers of $\nu\mathcal{H}$ but not necessarily to the zero section. So the intersection
of this foliation with the zero section will be a ``singular'' foliation. If the Haefliger structure is also transverse to the zero section, then it is called a {\it regular} Haefliger structure that comes from a genuine foliation. Two Haefliger structures are considered the same if they have the same normal bundle and if they coincide in a neighborhood of the zero section. 
\end{defn}
\begin{figure}[h]
  \centering
    \includegraphics[width=0.5\textwidth]{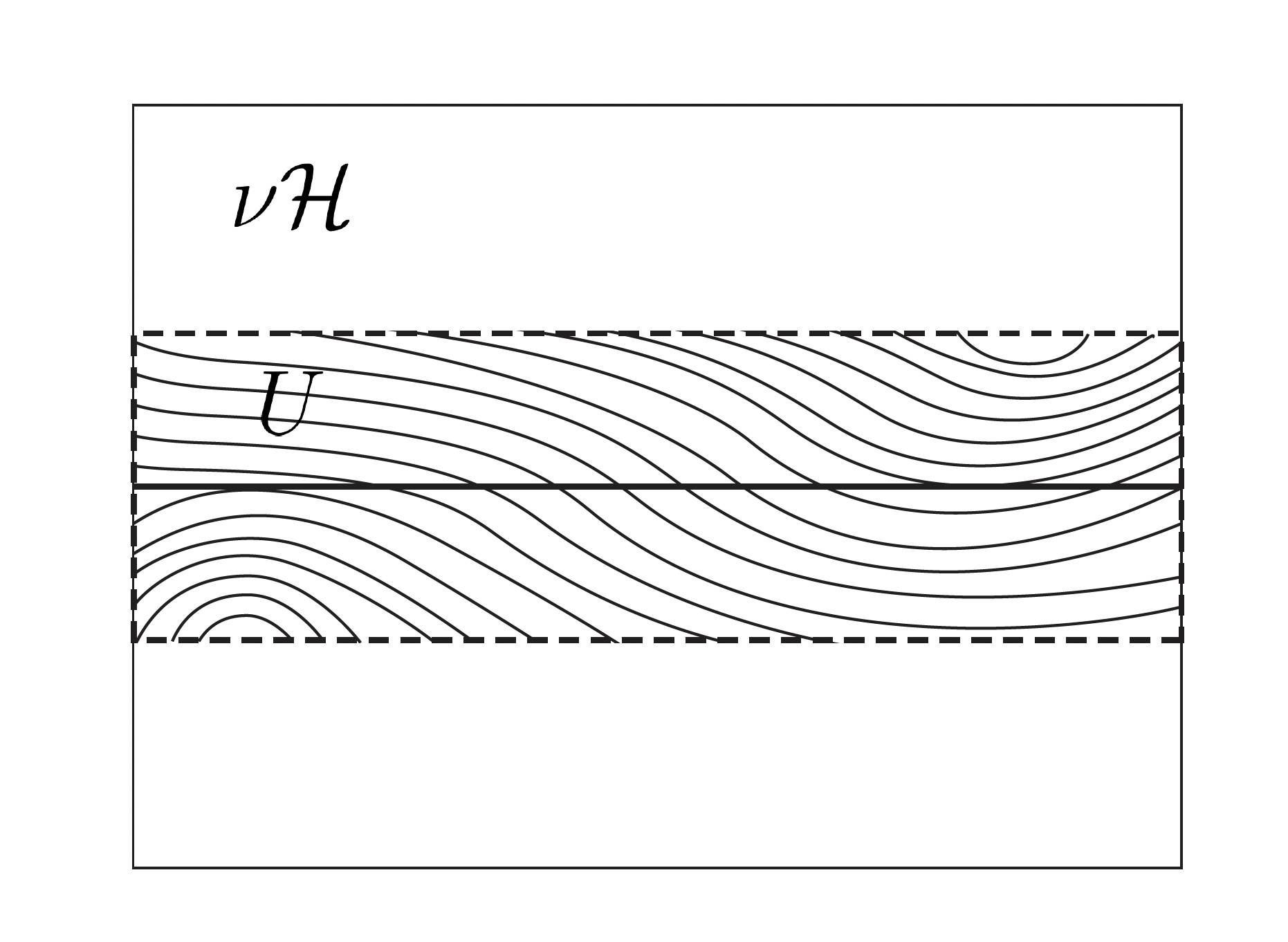}
    \caption{Haefliger structure}

\end{figure}

So one can think of a Haefliger structure as a germ of a foliated bundle near its zero section. We say two codimension $k$ Haefliger structures $\mathcal{H}_1$ and $\mathcal{H}_2$ on $M$ are {\it concordant}, if there is a codimension $k$ Haefliger structure $\mathcal{H}$ on $M\times [0,1]$ that restricts to $\mathcal{H}_i$ on $M\times \{i\}$ for $i=0,1$. A fundamental question is when a Haefliger structure is concordant with a regular one.   Using tools of homotopy theory, Haefliger constructed (\cite{MR0285027})  a classifying space $\mathrm{B}\Gamma^r_k$ for codimension $k$ $C^r$-Haefliger structures such that the concordance classes of codimension $k$ $C^r$-Haefliger structures on $M$ are in bijection with continuous maps  $M\to\mathrm{B}\Gamma^r_k$ up to homotopy. There is also a classifying space $\overline{\mathrm{B}\Gamma}^r_k$ that classifies Haefliger structures $\mathcal{H}$ with trivialized normal bundle $\nu\mathcal{H}$ (i.e. $\nu\mathcal{H}$ is a product bundle). 

Using Gromov's h-principle and a result of A. Phillips, Haefliger showed for an open manifold $M$, the existence of a codimension $k$ $C^r$ foliation on $M$ can be reduced to the homotopy theory question of the existence of a $C^r$ Haefliger structure $\mathcal{H}$ on $M$ such that its normal bundle $\nu\mathcal{H}$ can be embedded into the tangent bundle $TM$.  

The existence and classification of foliations on closed manifolds were much harder and there was a string of constructions by Reeb,  Lickorish, Novikov, Wood, Lawson, Durfee, and Tamura in particular cases (see \cite{MR0343289} and references therein). Even the case of finding codimension-one foliations on odd-dimensional spheres was stubborn for many years until Durfee and Tamura found constructions for all odd-dimensional spheres. Thurston proved (\cite{MR4554427}) dramatic general results for the existence and classification of foliations on closed manifolds. In particular, he showed that the concordance classes of foliations on a closed manifold $M$ are in bijection with homotopy classes of Haefliger structures $\mathcal{H}$ together with concordance classes of bundle monomorphisms $\nu\mathcal{H}\to TM$. 

These spectacular h-principle theorems reduced the subtle question of classifying foliations to the homotopy theory of Haefliger classifying spaces. But the homotopy type of the Haefliger spaces remains mysterious. The main theme of this article is the way Thurston studied the homotopy type of Haefliger classifying spaces in (\cite{MR0339267}) by generalizing a result of Mather (\cite{MR0356085}) to relate the homotopy groups of $\mathrm{B}\Gamma^r_k$   to the group homology of compactly supported $C^r$-diffeomorphism group $\Diff^r_c(\bR^n)$. 

The Mather-Thurston theorem is interesting in at least two ways. From a technical point of view, the original proof by Thurston used what was known about the homotopy type of $\mathrm{B}\Gamma^r_k$ at the time and the idea of {\it fragmentation} to prove an intrinsically ``compactly supported'' statement. Unlike other h-principle theorems, the case of the open manifolds for the Mather-Thurston theorem was very nontrivial and it was later proved by Segal and McDuff (\cite{MR0516216, MR0569248}) after they developed their group completion theorem and homology fibration techniques. From a philosophical point of view, it was used in the opposite way than h-principle theorems are normally used. Mather and Thurston applied this theorem to study the homotopy type of a more ``flexible'' object  $\mathrm{B}\Gamma^r_k$ using algebraic information (like the perfectness or the simplicity) of a more ``rigid'' object of the discrete group of diffeomorphisms. 

There are already many excellent surveys on different aspects of foliations and diffeomorphism groups including dynamical systems, differential topology, differential geometry, and non-commutative geometry. In particular, commentary on foliations in Thurston's collected works (\cite{MR4554427}) beautifully summarizes his many contributions to the topic. In this article, we are interested in different aspects of the relationship 
\[
\left\{\begin{tabular}{c}
       \text{Algebraic properties of}\\ \text{diffeomorphism groups}
 \end{tabular}  \right\}\leftrightarrow\left\{\begin{tabular}{c}
       \text{Homotopical properties of }\\ \text{foliations and their invariants}
 \end{tabular}  \right\}
\]
and its new consequences on invariants of flat bundles.  
\section{Haefliger-Thurston's conjecture}
Recall that the space $\mathrm{B}\Gamma^r_k$ classifies codimension $k$ $C^r$-Haefliger structures $\mathcal{H}$ and $\overline{\mathrm{B}\Gamma}^r_k$ classifies  codimension $k$ $C^r$-Haefliger structures whose  normal bundles $\nu\mathcal{H}$ are trivial. They sit in a fibration sequence
\[
\overline{\mathrm{B}\Gamma}^r_k\to \mathrm{B}\Gamma^r_k\xrightarrow{\nu} \mathrm{B}\text{GL}_k(\bR),
\]
where $\nu$ is the map that classifies the normal bundle to Haefliger structures. The space $\mathrm{B}\text{GL}_k(\bR)$ is the classifying space of the topological group $\text{GL}_k(\bR)$ that classifies $k$-dimensional vector bundles.
\begin{conj}[Haefliger-Thurston]\label{HT0}
The space $\overline{\mathrm{B}\Gamma}^r_k$ is $2k$-connected (i.e. $\pi_i(\overline{\mathrm{B}\Gamma}^r_k)=0$ for $i\leq 2k$). 
\end{conj}
Haefliger used a result of A. Phillips for open manifolds to prove that $\overline{\mathrm{B}\Gamma}^{r}_k$ is $k$-connected for all $r$.  Thurston  (\cite{MR0339267})  used his generalization of Mather's homology isomorphism on compact manifolds to show that  $\overline{\mathrm{B}\Gamma}^{\infty}_k$ is $(k+1)$-connected and shortly after, Mather (\cite[Section 7]{MR0356129}) proved the same statement for  $\overline{\mathrm{B}\Gamma}^{r}_k$ for all regularities $r$ except  $r=k+1$. Up to $k$-connectedness, the proof uses the ``flexibility'' of $\overline{\mathrm{B}\Gamma}^r_k$ but for higher homotopy groups, as we shall explain, the only tool known is to use the Mather-Thurston theorem and homological properties of diffeomorphism groups as discrete groups. 

\begin{rem} As a consequence of Mather's acyclicity result (\cite{MR0288777}), it is known that $\overline{\mathrm{B}\Gamma}^0_k$ is contractible. A remarkable theorem of Tsuboi (\cite{MR1014925}) also shows that $\overline{\mathrm{B}\Gamma}^1_k$ is contractible. However, for $r>1$,  Roussarie in codimension one and Thurston in all codimensions used Godbillon-Vey invariants to show that $\overline{\mathrm{B}\Gamma}^r_k$ is not $(2k+1)$-connected. This led to an extensive study of secondary invariants of foliations which is not the focus of this article. But let us just mention that the non-vanishing of the secondary invariants of foliations is the only known way to detect nontrivial homotopy groups of $\overline{\mathrm{B}\Gamma}^r_k$.
\end{rem}

Given Thurston's results on regularizing Haefliger structures, \Cref{HT0} implies that $n-k$ dimensional plane fields on an $n$-dimensional manifold $M$, where $k$ is at least $(n-1)/2$, are up to homotopy integrable to  $C^r$-foliations. The case of $k=1, n=3$ is known for all regularities except $r=2$ and it goes back to Mather's homology isomorphism. Using Thurston's results, it implies that all $2$-plane fields on closed $3$-manifolds are up to homotopy integrable to a $C^r$ foliation for $r\neq 2$. 

 Let's focus on the {\it smooth} case so we drop the regularity $r$ from the superscripts. Haefliger already proved that $\overline{\mathrm{B}\Gamma}_1$ is simply connected. So, to show that it is $2$-connected, by the Hurewicz theorem it is enough to show that $H_2(\overline{\mathrm{B}\Gamma}_1;\bZ)=0$. This was first proved by Mather as a consequence of a more general result for codimension-one foliations. Here we sketch an argument using geometric ideas that go back to Claude Roger (\cite{MR0348767}) and Thurston. Later Meigniez generalized this geometric point of view to higher codimensions in \cite{MR4251434}. This geometric perspective unlike the other proofs, is more useful for the constructive statements in \Cref{CMT}.

\subsection{Finding a normal form for codimension-one Haefliger structures} For a topological group $G$, we let $G^{\delta}$ denote the same group with the discrete topology. Recall that the classifying space $\mathrm{B}G$ classifies the principal $G$-bundles and the classifying space $\mathrm{B}G^{\delta}$ classifies flat principal $G$-bundles. One way to think about flat bundles is that there exists a foliation on the total space that is transverse to the fibers. The space $\mathrm{B}G^{\delta}$ has $G^{\delta}$ as the fundamental group and has vanishing higher homotopy groups. So a continuous map $g\colon M\to \mathrm{B}G^{\delta}$ induces a representation $\rho\colon \pi_1(M)\to G^{\delta}$. Now let $\widetilde{M}$ denote the universal cover of $M$ and consider the horizontal foliation on $\widetilde{M}\times G$ whose leaves are $\widetilde{M}\times \{x\}$ for $x\in G$. Since this foliation is invariant under the diagonal action of $\pi_1(M)$ on  $\widetilde{M}\times G$, it induces a foliation on the quotient $\widetilde{M}\times_{\pi_1(M)} G$. So  we obtain a principal $G$ bundle
\[
G\to \widetilde{M}\times_{\pi_1(M)} G\to M,
\]
with a foliation on the total space which is transverse to the fibers. So in this sense, $\mathrm{B}G^{\delta}$ classifies flat $G$-bundles, and also the group (co)homology of the group $G^{\delta}$ is the same as the (co)homology of the classifying space $\mathrm{B}G^{\delta}$. We can think of the cohomology of the classifying space $\mathrm{B}G^{\delta}$ as the ring of characteristic classes of flat $G$-bundles which is another reason to study the cohomology of $\mathrm{B}G^{\delta}$. We shall mention some non-vanishing results about characteristic classes of flat $G$-bundles in \Cref{eqMT} but here we are mostly concerned about the relation between the group homology of diffeomorphism groups and the homotopy groups of Haefliger classifying spaces.

To prove that $H_2(\overline{\mathrm{B}\Gamma}_1;\bZ)=0$, first we construct a map
\begin{align}\label{H_2}
H_1(\BcdDiff \bR;\bZ)&\xrightarrow{}H_2(\overline{\mathrm{B}\Gamma}_1;\bZ),
\end{align}
which turns out to be an isomorphism.
To define this so-called ``suspension'' map, first note that the group $\cdDiff \bR$ is isomorphic to its subgroup $\cdDiff {(0,1)}$ and it is easy to see that the inclusion induces a group homology isomorphism. So it is enough to define the map $H_1(\BcdDiff {(0,1)};\bZ)\to H_2(\overline{\mathrm{B}\Gamma}_1;\bZ)$. Let $I$ be the closed unit interval and consider the horizontal foliation on $I\times I$. Given $f\in \cdDiff {(0,1)}$, the horizontal foliation on $I\times I$ induces a foliation $\mathcal{H}_f$ on
\[
S^1\times I=\frac{I\times I}{\{0\}\times x \sim \{1\}\times f(x)}\cdot
\]
The foliation $\mathcal{H}_f$ is horizontal near the two boundary components. So by contracting the two circle boundary components, we obtain a singular foliation (aka Haefliger structure) $\hat{\mathcal{H}}_f$ on $S^2$ (see  Figure 5) with two center singularities. 
\begin{figure}[h]\label{suspension}
  \centering
    \includegraphics[width=0.4\textwidth]{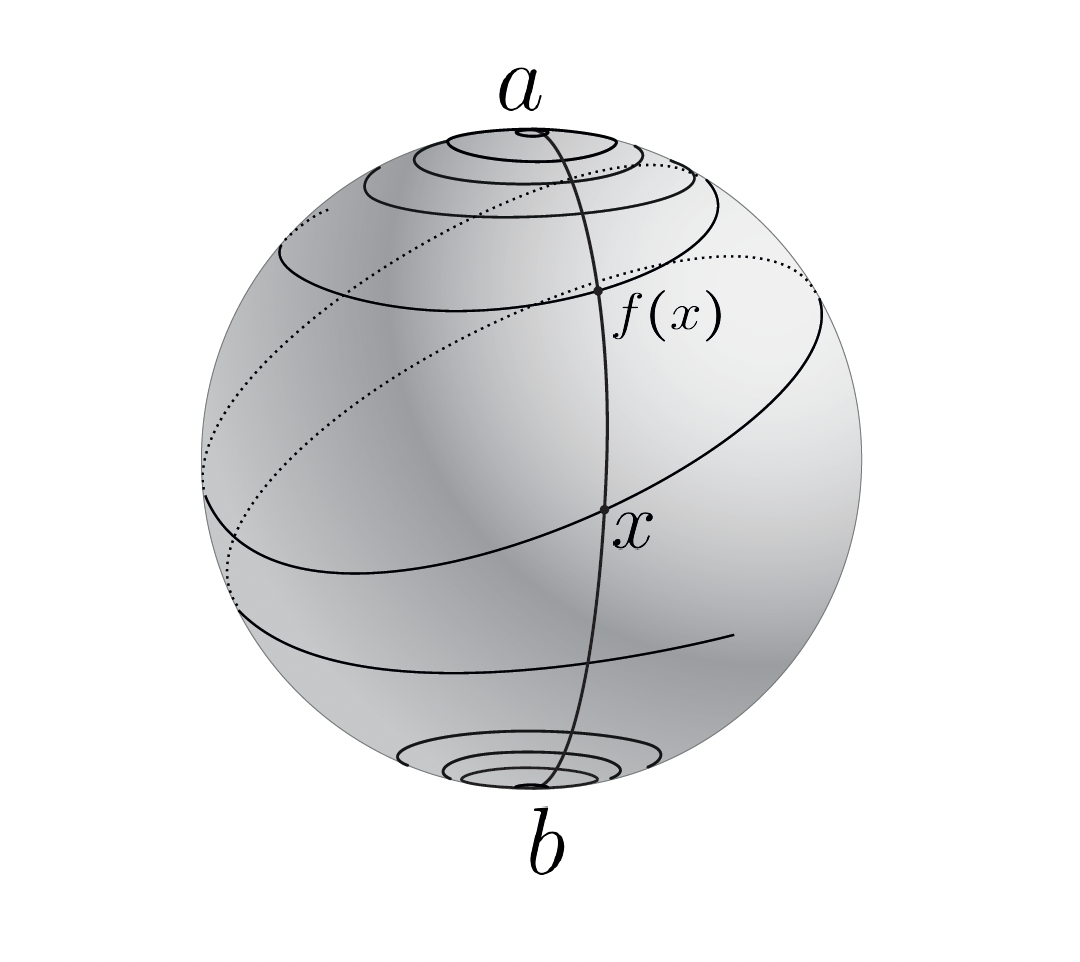}
    \caption{The Haefliger structure $\hat{\mathcal{H}}_f$ obtained by the suspension map}

\end{figure}
So the surjection of the suspension map \ref{H_2} is equivalent to the following.
\begin{thm}[Mather]\label{normalform} Each codimension-one Haefliger structure on $S^2$ is concordant with a Haefliger structure $\hat{\mathcal{H}}_f$ given by the suspension of a diffeomorphism $f\in \Diff_c {((0,1))}$. 
\end{thm}
 Note that an oriented codimension-one Haefliger structure  on $S^2$ is induced by a map $S^2\to \overline{\mathrm{B}\Gamma}_1$ so it is an element  in $\pi_2(\overline{\mathrm{B}\Gamma}_1)=H_2(\overline{\mathrm{B}\Gamma}_1;\bZ)$. This element can be represented by a germ of a foliation $\mathcal{F}$ on $S^2\times \bR$ near the zero section $S^2\times \{0\}$. On each foliation chart $U_i$, we have a submersion $f_i\colon U_i\to \bR$ (i.e. $Df_i(x)\neq 0$ for all $x\in U_i$) such that $\mathcal{F}|_{U_i}$ is given by the pullback of point foliation on $\bR$ via $f_i$. Since $\mathcal{F}$ is not necessarily transverse to the zero section $S^2\times\{0\}$, the restrictions of functions $f_i$ to the zero section have singularities. One can change $\mathcal{F}$ in its concordance class to make all the singularities Morse singularities. So we end up having either center or saddle singularities.

\begin{figure}[h]\label{singularity}
  \centering
    \includegraphics[width=0.5\textwidth]{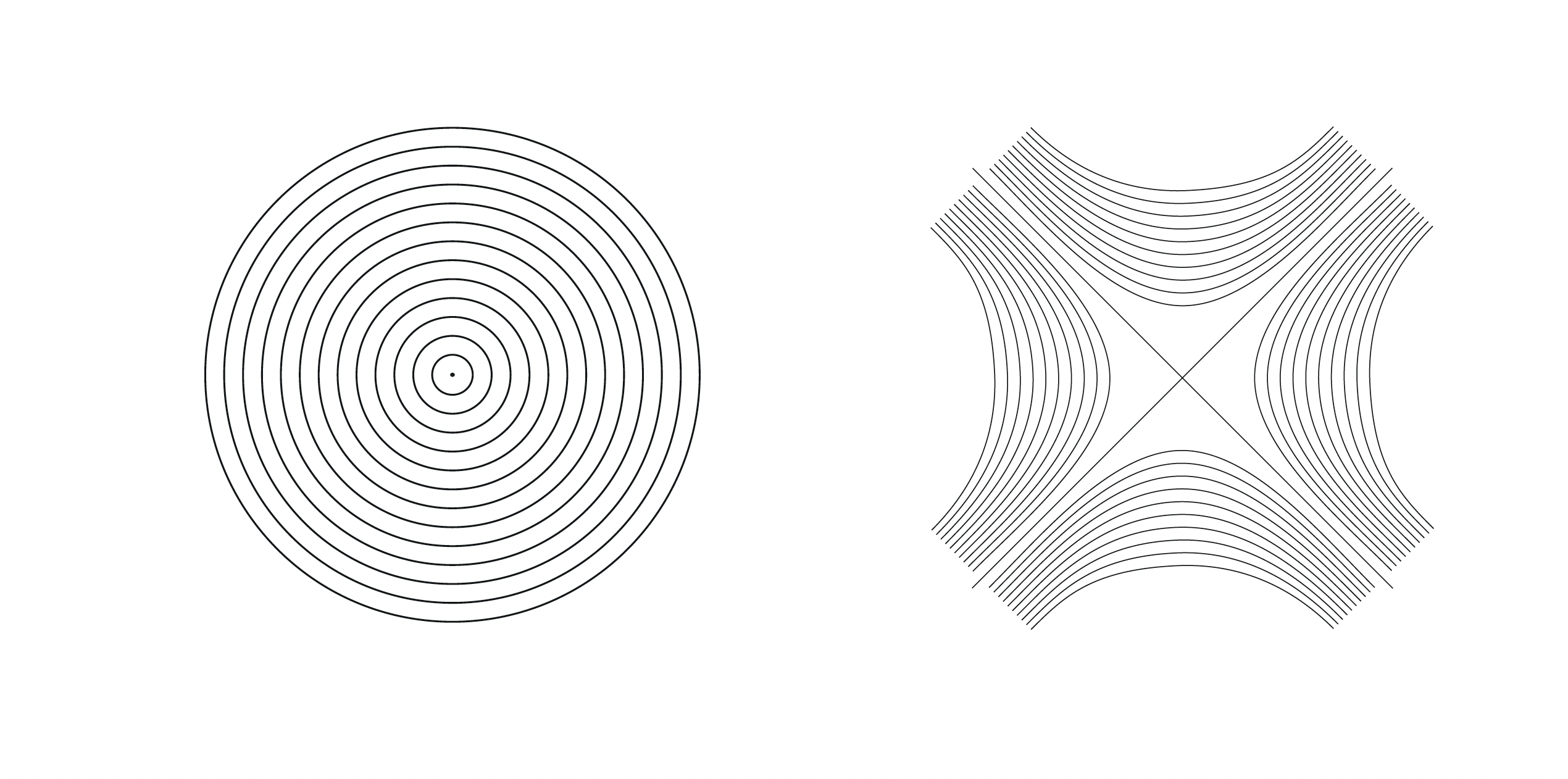}
    \caption{Center and saddle singularities}

\end{figure}

Let $\mathcal{H}$ be a Haefliger structure on $S^2$ with a trivialized normal bundle that has only Morse singularities. Then there is a standard strategy to remove singularities by pairing each saddle singular point to a center. Claude Roger observed that one can choose a path between a saddle and a center and change the Haefliger structure up to concordance only in a neighborhood of the path to remove both singularities. To prove \Cref{normalform}, we continue this process until we are left with two center singularities on $S^2$, which is guaranteed by the Hopf index theorem. Such a Haefliger structure as in Figure 5 is induced by the suspension of a diffeomorphism. 


\subsection{Bordism of Haefliger structures and the simplicity trick} Given the normal form in \Cref{normalform} for codimension-one Haefliger structures on $S^2$, to prove $H_2(\overline{\mathrm{B}\Gamma}_1;\bZ)=0$, it is enough to show that for each $f\in \Diff_c {((0,1))}$ the Haefliger structure  $\hat{\mathcal{H}}_f$ on $S^2$ can be extended to a $3$-manifold bounding $S^2$. Surprisingly the only proof known for this geometric question is to use a deep result of Herman that the group of orientation preserving diffeomorphisms of the circle $\Diff_0{(S^1)}$ is a perfect group meaning that it is equal to its commutator subgroup. This is against the philosophy of h-principle though, since one might expect that it ought to be easier to compute $H_2(\overline{\mathrm{B}\Gamma}_1;\bZ)$ than to understand the more rigid object of diffeomorphism groups!

Given Herman's theorem, if we show that  $\hat{\mathcal{H}}_f$ is cobordant to zero (i.e. $\hat{\mathcal{H}}_f$ on $S^2$ can be extended to a $3$-manifold bounding $S^2$) for all $f\in \Diff_c {((0,1))}$, we conclude that $\overline{\mathrm{B}\Gamma}_1$ is $2$-connected. Roger used the perfectness of $\Diff_0{(S^1)}$ to show that the null-bordism for $\hat{\mathcal{H}}_f$ exists. But let us mention a trick due to Thurston that builds a specific null-bordism. Choosing an ``economical'' bordism between Haefliger structures will show up in the last section again. 

The perfectness of $\Diff_0(S^1)$  implies that it is a simple group by Epstein's theorem (\cite{MR0267589}), meaning it has no nontrivial normal subgroup.  Thurston's {\it trick} uses the simplicity of $\Diff_0{(S^1)}$ to show that $\hat{\mathcal{H}}_f$ is cobordant to zero. The first step is to show that $\hat{\mathcal{H}}_f$ is cobordant to a foliation on the torus $S^1\times S^1$ that is transverse to the circles $\{x\}\times S^1$ as follows.

\begin{figure}[h]\label{cobordism}
  \centering
    \includegraphics[width=0.5\textwidth]{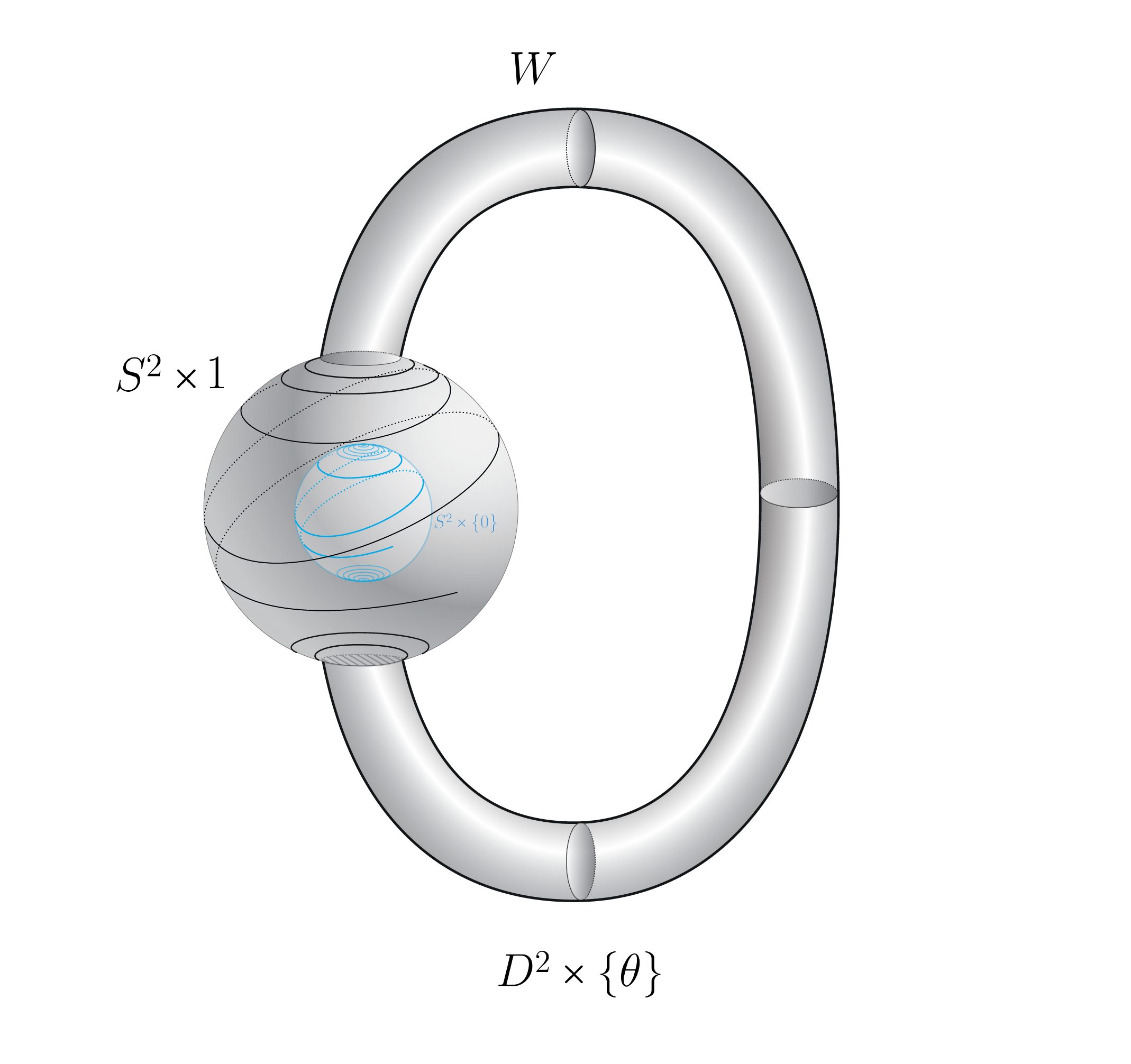}
    \caption{Bordism of Haefliger structures}
\end{figure}
 Let $W$ be a $3$-manifold obtained by gluing $D^2\times I$ to $S^2\times I$ as in Figure 7. Thus $\partial W$ is diffeomorphic to $S^2\coprod S^1\times S^1$. A Haefliger structure on $W$ is given by the product Haefliger structure $\hat{\mathcal{H}}_f\times I$ on $S^2\times I$ and the codimension one foliation given by leaves $D^2\times \{x\}$ on $D^2\times I$. Note that the boundary of $W$ is the union of a sphere and a torus. This gives a bordism between $\hat{\mathcal{H}}_f$ on its sphere boundary component and a foliation on the torus boundary component that is transverse to $\{x\}\times S^1$ for all $x\in S^1$. Hence, this foliation $\mathcal{F}_{\phi}$ on the torus is induced by suspending the horizontal foliation on the cylinder $S^1\times I$ by a diffeomorphism $\phi\in \Diff_0{(S^1)}$.

The idea of Thurston is to consider the subset of the circle diffeomorphisms $\phi$ for which the suspension $\mathcal{F}_{\phi}$ can be extended to a foliation on a solid torus. Then it is not hard to see that this subset is in fact a {\it normal subrgroup} of $\Diff_0{(S^1)}$. Given that $\Diff_0{(S^1)}$ is a simple group, it is enough to show that this contains at least one non-trivial diffeomorphism. Suppose that $\phi$ is a nontrivial rotation of the circle. Thurston extends $\mathcal{F}_{\phi}$ into a solid torus by standard construction of {\it spinning leaves} around an embedded torus inside the solid torus and filling in the torus leaf by a Reeb foliation as in Figure 8. We conclude that $\hat{\mathcal{H}}_f$ is cobordant to a foliation on the torus that is nullcobordant; therefore $\hat{\mathcal{H}}_f$ is indeed nullcobordant. 

\begin{figure}[h]\label{Reeb}
  \centering
    \includegraphics[scale=0.25]{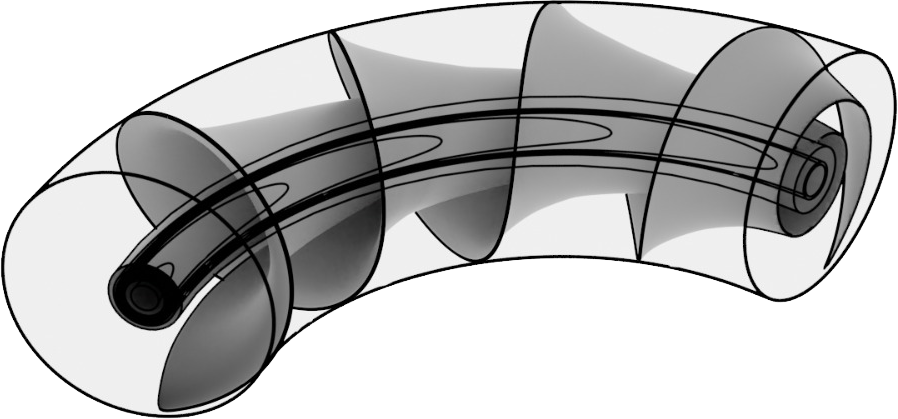}
    \caption{Part of a solid torus with spinning leaves and a Reeb component}
\end{figure}
As we shall see in the next section, it is a consequence of a more general theorem of Mather that $\Diff^r_0{(S^1)}$ is simple for $r\neq 2$. So the above construction also for $r\neq 2$ shows that $C^r$ codimension-one Haefliger structures on $S^2$ are cobordant to the trivial Haefliger structure. 
\begin{problem}
Is $\Diff^2_0{(S^1)}$ a simple group?
\end{problem}
By Mather's theorem \ref{M} below, If one could find a geometric construction without using perfectness or simplicity of diffeomorphism groups to show that suspension Haefliger structures $\hat{\mathcal{H}}_f$ for all $f\in \Diff_c^r {((0,1))}$ is null-cobordant, it would give a geometric proof of the simplicity of $\Diff_0^r{(S^1)}$ in all regularities. But so far, no such construction is known without using the algebraic properties of these groups, which is against the h-principle philosophy!

Mather proved a more general theorem for codimension-one foliations that also gives a non-geometric proof that $\overline{\mathrm{B}\Gamma}_1$ is $2$-connected. Consider the classifying space $\BcdDiff \bR$ and the trivial bundle $\BcdDiff \bR\times \bR\to \BcdDiff \bR$. This is a ``universal compactly supported foliated $\bR$-bundle''  in the sense that for each foliated $\bR$-bundle $M\times \bR\to M$ over a manifold $M$ where the foliation on $M\times \bR$ is horizontal outside of a compact subset, is a pullback of this universal bundle via a map $M\to \BcdDiff \bR$. Hence, there is classifying map $\BcdDiff \bR\times \bR\to \overline{\mathrm{B}\Gamma}_1$. The adjoint of this map \footnote{Recall that the adjoint of a continuous map $X\times Y\to Z$ is a continuous map $X\to \text{Map}(Y,Z)$.} is a map  $\BcdDiff \bR\to \Omega \overline{\mathrm{B}\Gamma}_1$ where $\Omega \overline{\mathrm{B}\Gamma}_1$ is the  based loop space of $ \overline{\mathrm{B}\Gamma}_1$.
\begin{thm}[Mather]\label{M} The adjoint map
\[
\BcdrDiff \bR\to \Omega \overline{\mathrm{B}\Gamma}^r_1
\]
induces a homology isomorphism.
\end{thm}
In particular, Mather's theorem implies that the suspension map (\ref{H_2})
\[
H_1(\BcdDiff \bR;\bZ)\to H_1(\Omega \overline{\mathrm{B}\Gamma}_1;\bZ)\cong H_2(\overline{\mathrm{B}\Gamma}_1;\bZ),
\]
is in fact an isomorphism. Mather's theorem also has a version for $\drDiffo {S^1}$ which implies that there is an isomorphism
\[
H_1(\BcdrDiff \bR;\bZ)\xrightarrow{\cong}H_1(\BdrDiffo {S^1} ;\bZ).
\]
Since $H_1(\mathrm{B}G^{\delta};\bZ)$ is isomorphic to the abelianization $G/[G,G]$, Herman's perfectness theorem implies that $H_1(\BdDiffo {S^1};\bZ)=0$. Therefore, Mather's homology isomorphism implies that $H_2(\overline{\mathrm{B}\Gamma}_1;\bZ)=0$.
\section{The Mather-Thurston theorem} Mather's original proof of \Cref{M} is quite involved. He showed that the delooping of $\BcdDiff \bR$ using  Quillen's delooping techniques in Algebraic K-theory is homotopy equivalent to $\overline{\mathrm{B}\Gamma}_1$. This step was not easy to generalize to higher dimensions. But Thurston found a generalization using a method of {\it fragmentation} for compactly supported diffeomorphisms. We shall explain this relation between  {\it fragmentation} and the desired delooping in \Cref{Proof:MT}.

To state his theorem, let $M$ be an $n$-dimensional smooth manifold and let $\overline{\BDiff_c^r(M)}$ be the homotopy fiber of the map
\[
\eta\colon \BcdrDiff M\to \BDiff^r_c(M),
\]
which classifies trivialized $M$-bundles equipped with foliations on the total space transverse to the fibers such that the foliation is {\it compactly supported} i.e. it is horizontal outside of some compact subset of the total space. Let $\Diff_{c,0}^r(M)$ be the identity component of $\Diff_{c}^r(M)$ and let $\widetilde{\Diff}_{c,0}^r(M)$ be its universal cover. It turns out that the fundamental group of $\overline{\BDiff_c^r(M)}$ is isomorphic to $\widetilde{\Diff}_{c,0}^r(M)^{\delta}$. Now the natural map $\widetilde{\Diff}_{c,0}^r(M)^{\delta}\to \Diff_{c,0}^r(M)$ gives a representation of $\pi_1(\overline{\BDiff_c^r(M)})$ into $\Diff_{c,0}^r(M)$ that induces a trivialized $M$-bundle $\overline{\BDiff_c^r(M)}\times M\to \overline{\BDiff_c^r(M)}$; we call it the universal foliated trivialized $M$-bundle. This universal foliation on the total space induces a diagram commutative 
 up to homotopy
\[
 \begin{tikzpicture}[node distance=1.8cm, auto]
  \node (A) {$\overline{\BDiff_c^r(M)}\times M$};
  \node (B) [right of=A, node distance=3cm] {$\mathrm{B}\text{GL}_n(\bR),$};
  \node (C) [above of= B ] {$\mathrm{B}\Gamma_n^r$};  
   \draw [->] (A) to node {$\theta$}(B);
  \draw [->] (C) to node {$\nu^r$}(B);
  \draw [->] (A) to node {$\alpha$}(C);
\end{tikzpicture}
\]
where $\alpha$ classifies the codimension $n$ foliation on the total space. The map $\theta$ classifies the normal bundle to the foliation and since the foliation is transverse to the fibers, this map is induced by the map $\tau_M\colon M\to  \mathrm{B}\text{GL}_n(\bR)$ which classifies the tangent bundle of $M$.

The adjoint of the map $\alpha$ induces a map from $\overline{\BDiff_c^r(M)}\times M$ to the space of lifts of $\tau_M$ to $\mathrm{B}\Gamma_n^r$. One can think about this space of lifts as follows. Let $\tau_M^*(\nu^r)$ be the bundle induced by pulling back the fibration $\nu^r$ via $\tau_M$. Fix a section $s_0$ of this bundle. For any other section $s$, we can define the support of $s$ to be the closure of the set of points $x\in M$ where $s(x)\neq s_0(x)$. Let $\text{Sect}_c(\tau_M^*(\nu^r))$ be the subspace of compactly supported sections of the bundle $\tau_M^*(\nu^r)\to M$. 

   
Therefore, the adjoint of $\alpha$ gives a map
\[
f\colon\overline{\BDiff_c^r(M)}\to \text{Sect}_c(\tau_M^*(\nu^r)).
\]
\begin{thm}[Mather-Thurston]\label{MT}
The map $f$ is an acyclic map; in particular, it induces a homology isomorphism with integer coefficients.
\end{thm}
A {\it non-compactly supported} version of this theorem also holds. But surprisingly, even the case of an open ball was proved later by McDuff and Segal. As we shall explain, Thurston's fragmentation method only works for the compactly supported version.

 Thurston gave a few applications of this theorem. To mention two of them, note that the target of the map $f$ above, which is a section space of a fiber bundle is more amenable to homotopy theory. The fiber of this bundle is $\overline{\mathrm{B}\Gamma}^{r}_n$, which was known by Haefliger to be $n$-connected. Using obstruction theory, one can see that $$H_1(\text{Sect}_c(\tau_M^*(\nu^r));\bZ)=\pi_{n+1}(\overline{\mathrm{B}\Gamma}^{r}_n),$$ and in particular it does not depend on $M$. Therefore, $H_1(\overline{\BDiff_c^r(M)};\bZ)$  does not depend on $M$ either. Recall that  $\pi_1(\overline{\BDiff_c^r(M)})$ is isomorphic to $\widetilde{\Diff}_{c,0}^r(M)^{\delta}$ which was known by a deep result of Herman and Mather to be perfect for the case of the torus $M=T^n$ and $r\neq n+1$. Therefore, the universal cover of $\Diff_{c,0}^r(M)$ is perfect for all $M$ and $r\neq n+1$ which implies the same for $\Diff_{c,0}^r(M)$ itself. As a result of Epstein, the perfectness for $\Diff_{c,0}^r(M)$ implies that these groups are simple. 
 
 Another application related to the theme of this article is that the perfectness result implies that $\pi_{n+1}(\overline{\mathrm{B}\Gamma}^{r}_n)=0$, so given Haefliger's result that $\overline{\mathrm{B}\Gamma}^{r}_n$ is $n$-connected, the connectivity of $\overline{\mathrm{B}\Gamma}^{r}_n$ is improved by one for $r\neq n+1$.
 
 Using the Mather-Thurston theorem, one can show that Haefliger-Thurston's conjecture is equivalent to the following bundle theoretic statement.
 \begin{conj}[Haefliger-Thurston]The map
 \[
\eta\colon \BcdrDiff M\to \BDiff^r_c(M),
\]
induces an isomorphism on homology in degrees less than $n+1$  and is a surjection on homology in degree $n+1$. 
 \end{conj}
 \begin{rem}
 In low regularities when $r=0$ using the acyclicity result of Mather and for $r=1$ using the acyclicity result of Tsuboi, we know that $\eta$  in fact induces a homology isomorphism in all degrees. But for $r>1$, the homological degrees in the conjecture cannot be improved.
 \end{rem}
 
 Geometrically, this conjecture is equivalent to the following. For every $C^r$ $M$-bundle $M\to E\to B$ where $B$ is a closed manifold and $\text{dim}(B)\leq \text{dim}(M)$, there exists a bordism  $W$ from $B$ to another manifold $B'$ and a $C^r$ $M$-bundle $M\to K\to W$ such that when it is restricted to $B$, it is isomorphic to $E\to B$ and when it is restricted to $B'$, it is a flat $M$-bundle i.e. it is induced by a representation $\pi_1(B')\to \Diff_c^{r}(M)^{\delta}$. 
 
 For non-compactly supported diffeomorphisms of open manifolds, the range of isomorphism in the conjecture is expected to be better. For example, for $M=\bR^n$, Segal proved that the natural map
 \[
 \BdrDiff {\bR^n}\to \mathrm{B}\Gamma^r_n,
 \]
 induces a homology isomorphism. It is easy to see that Segal's theorem implies that the map 
 \[
  \BdrDiff {\bR^n}\to \BDiff^r(\bR^n),
 \]
 induces a homology isomorphism already in degrees less than $n+1$. Linearizing diffeomorphisms implies that $\Diff^r(\bR^n)\simeq \text{GL}_n(\bR)$ for $r>0$, so as a consequence, we know the homology of $\BdrDiff {\bR^n}$ up to degree $n$. In particular, Segal's theorem has the following corollary.
 
 \begin{cor} 
 The group $\drDiffo {\bR^n}$ is perfect for all $r$.
 \end{cor} 
 This proof is in line with the h-principle philosophy and works for all regularities. McDuff and Schweitzer later found an algebraic proof for the perfectness of $\drDiffo M$ using ``Ling's factorization'' when $M$ is the interior of a compact manifold with a boundary. 
 
McDuff also proved the analog of the Mather-Thurston theorem for volume-preserving diffeomorphisms and Segal's theorem for the volume-preserving case for $n>2$. Curiously, the case and $n=2$ is still open.
 \begin{problem}Let  $\drDiff {\bR^2, \omega}$ be the group of volume preserving diffeomorphisms of $\bR^2$ with the standard volume form $\omega$. Is $\drDiff {\bR^2, \omega}$ perfect? Is the natural map from $\BdrDiff {\bR^2, \omega}$ to the corresponding Haefliger space $ \mathrm{B}\Gamma^{\omega}_2$ a homology isomorphism?
 \end{problem}
The case of volume-preserving homeomorphisms is also very interesting.  There has been a recent breakthrough in proving the failure of the perfectness of volume-preserving homeomorphisms of the $2$-disk by Daniel Cristofaro-Gardiner, Vincent Humilière, and Sobhan Seyfaddini, which is beyond the scope of this article. 
\subsection*{The curious case of Piecewise Linear homeomorphisms}The subgroups of piecewise linear (PL for short) homeomorphisms of the line have been a rich source for interesting finitely generated groups with surprising algebraic and dynamical properties. However, not much is known about the algebraic properties of the PL homeomorphisms of higher dimensional PL manifolds. Although they are more combinatorial in nature, the analytical tools for diffeomorphism groups and the Mather infinite repetition trick for homeomorphisms (\cite{MR0288777}) are not available for PL homeomorphisms. So the following basic question due to Epstein (\cite{MR0267589}) is still open.
\begin{problem}(Epstein) Let $M$ be a PL manifold. Is $\text{PL}_0(M)$, the group of PL homeomorphisms of $M$ that are isotopic to the identity, a simple group? 
\end{problem}
By Epstein's result, perfectness and simplicity are equivalent in this case and he proved that $\text{PL}_c(\bR)$ and $\text{PL}_0(S^1)$ are perfect by observing that in dimension one,  PL homeomorphisms are generated by certain ``typical elements'' and those typical elements can be easily written as commutators. To generalize his argument to higher dimensions, he suggested the following approach.

Let $B$ be a ball in $\bR^n$. It is PL homeomorphic to $S^{n-2}\star [0,1]$, the join of $S^{n-2}$ with $[0,1]$. Note that for PL manifolds $M$ and $N$, a PL homeomorphism of $N$ extends naturally to a PL homeomorphism of the join $M\star N$.  A {\it glide} homeomorphism of the ball $B$ is a PL homeomorphism that is induced by the extension of a compactly supported PL homeomorphism of $(0,1)$ to a PL homeomorphism of $S^{n-2}\star [0,1]$. For a PL $n$-manifold $M$, a  glide homeomorphism $h\colon M\to M$ is the extension by the identity of a glide homeomorphism supported in a PL embedded ball $B\hookrightarrow M$. 
However, it is not known if ${\textnormal{\text{PL}}}_{c,0}(M)$ is generated by glide homeomorphisms in all dimensions.

On the other hand, Greenberg started the program to study the Haefliger classifying space for PL foliations and he described them inductively in terms of the codimension of the foliation and observed interesting connections to algebraic $K$-theory of real numbers.  Greenberg's works (see \cite{MR1200422} and references therein)  describe an inductive method to build $\overline{\mathrm{B}\Gamma}_n^{\textnormal{\text{PL}}}$ from classifying spaces of matrix groups made discrete. In codimension one, he shows that $\overline{\mathrm{B}\Gamma}_1^{\textnormal{\text{PL}}}$ is homotopy equivalent to $K(\bR_+,1)*K(\bR_+,1)$ the join of two Eilenberg-MacLane spaces where $\bR_+$ is multiplicative group of the positive reals. In higher codimensions, he describes a complicated diagram whose homotopy colimit is homotopy equivalent to $\overline{\mathrm{B}\Gamma}_n^{\textnormal{\text{PL}}}$. Each space in Greenberg's diagram is related to classifying spaces of PL foliations of lower codimensions or classifying spaces of certain matrix groups made discrete. Already in codimension $2$, his model is difficult to do calculations with; in particular one cannot read off from the lower homotopy groups of $\overline{\mathrm{B}\Gamma}_2^{\textnormal{\text{PL}}}$. I used Suslin's work on the algebraic $K$-theory of $\bR$ in \cite{Nariman} to prove the analog of Haefliger-Thurston's conjecture for codimension $2$ PL foliations.
\begin{thm}(Nariman)
The classifying space $\overline{\mathrm{B}\Gamma}_2^{\textnormal{\text{PL}}}$ is  $4$-connected.
\end{thm} 

Greenberg (\cite{MR1200422}) and Gelfand-Fuks independently conjectured that the analog of the Mather-Thurston theorem holds for PL foliations. I used Thurston's fragmentation method in \cite{MR4580303} to prove the Mather-Thurston theorem in the PL case.  Using this version of Mather-Thurston for PL foliations of codimension $2$ and the fact that $\pi_3(\overline{\mathrm{B}\Gamma}_2^{\textnormal{\text{PL}}})=0$ solves  Epstein's question about PL homeomorphisms in dimension $2$. 

 \begin{thm}(Nariman)
Let $\Sigma$ be an oriented compact surface, possibly with a boundary. The identity component ${\textnormal{\text{PL}}}_0(\Sigma,\text{rel }\partial)$ is a simple group as a discrete group.
\end{thm}
This is a homotopy theoretical proof of the simplicity of ${\textnormal{\text{PL}}}_0(\Sigma,\text{rel }\partial)$ which is more in line with the h-principle philosophy and it is unlike the smooth case. Still, it would be interesting to prove the perfectness and simplicity of ${\textnormal{\text{PL}}}_0(\Sigma,\text{rel }\partial)$  group-theoretically.

\section{Immersion theoretic method vs Fragmentation method}\label{Proof:MT} Thurston apparently had three different proofs for \Cref{MT}. He gave lectures at Harvard on one of his proofs and Mather worked out the details and wrote it down (see \cite{MR4580303} and references therein). He said in an email to the author ``The lectures were sketchy and it was really hard to write up the proof. I spent 14 months on it''. McDuff and Segal  used their group completion theorem and homology fibration techniques to give another proof in the spirit of immersion theory. Their approach has the advantage that it also works for open manifolds. 

Recently Meigniez (\cite{MR4251434}) found a geometric proof similar to the proof of the isomorphism \ref{H_2} by Roger in which he changes Haefliger structure up to concordance to obtain a ``normal'' form with specified models for their singularities and he then does surgery on singularities to make the Haefliger structure cobordant to a foliated bundle (see Figure 7). This geometric approach is more useful for constructive statements like the ones in \Cref{CMT}. 

We shall sketch the main idea of how Thurston proved that the natural map
\begin{equation}\label{T}
\overline{\BDiff_c^r(\bR^n)}\to \Omega^n \overline{\mathrm{B}\Gamma}_n^r,
\end{equation}
where $\Omega^n \overline{\mathrm{B}\Gamma}_n^r$ is an $n$-fold loop space, is a homology isomorphism. This is in fact the main step in proving a compactly supported version of the Mather-Thurston theorem for all manifolds.

Let us first recall the immersion theoretic method that is used in Smale-Hirsch's theorem. Suppose $F_1, F_2\colon (\mathsf{Mfld}^{\partial}_n)^{op}\to \mathsf{S}$ are functors from the category of smooth $n$-manifolds (possibly with nonempty boundary) with smooth embeddings as morphisms to a convenient category of spaces $\mathsf{S}$.  To prove a natural transformation $T\colon F_1\to F_2$  induces a weak homotopy equivalence (i.e. it induces isomorphisms on all homotopy groups) $F_1(N)\simeq F_2(N)$ for all $N\in \mathsf{Mfld}^{\partial}_n$, we fix a handle decomposition of $N$ and induct on the number of handles as follows. For example, in the Smale-Hirsch theorem there is an ambient manifold $L$ whose dimension is larger than $n$, the functor $F_1(N):= \text{Imm}(N, L)$ is the space of immersions from $N$ to $L$, $F_2(N):=\text{Mon}(TN, TL)$ is the space of bundle monomorphisms between tangent bundles, and the natural transformation $T$ is induced by the differential of immersions.

First, we show that $F_1(\bR^n)\to F_2(\bR^n)$ is a weak equivalence. This is normally an easy step by ``zooming in'' arguments. Now suppose we know that $F_1(M)\to F_2(M)$ is a weak equivalence and  $N$ is a manifold obtained from $M$ by attaching a $k$-handle $D^k\times D^{n-k}$. So $M\cap (D^k\times D^{n-k})$ is $A^k\times D^{n-k}$ where $A^k$ is a neighborhood of $S^{k-1}$ in $D^k$. In the case of Smale-Hirsch's theorem for $i=1, 2$, we have pull-back squares
\begin{equation}\label{d1}
\begin{gathered}
\begin{tikzpicture}[node distance=3cm, auto]
  \node (A) {$F_i(N)$};
  \node (B) [right of=A] {$F_i(M)$};
  \node (C) [below of=A, node distance=1.2cm] {$ F_i(D^k\times D^{n-k})$};  
  \node (D) [below of=B, node distance=1.2cm] {$F_i(A^k\times D^{n-k}),$};
  \draw[->] (C) to node {$$} (D);
  \draw [->] (A) to node {}(C);
  \draw [->] (A) to node {$$} (B);
  \draw [->] (B) to node {$$} (D);
\end{tikzpicture}
\end{gathered}
\end{equation}
where the maps are given by restrictions. If the horizontal maps are Serre fibrations then the squares become {\it homotopy} pull-back squares. Therefore, by induction knowing the weak equivalence for $M$, $D^k\times D^{n-k}$ and $A^{k}\times D^{n-k}$, we can compare the two homotopy pull-back squares and conclude that $F_1(N)\to F_2(N)$ is a weak-equivalence. 

Normally in h-principle theorems, the functor $F_2$ has nice homotopical properties, for example it is a mapping space or a section space of a bundle over the manifold. In these cases, it is easy to see that the restriction maps similar to the diagram \ref{d1} are Serre fibrations. And the hard step is to prove the same for a more geometric functor $F_1$. In the case of Mather-Thurston, the base case was nontrivial and it was proved by Segal that $F_1(\bR^n)\to F_2(\bR^n)$ is a homology isomorphism. McDuff showed that the restriction maps for $F_1$ in the Mather-Thurston theorem are {\it homology} fibrations.   

Thurston proved the compactly supported version without knowing Segal's theorem for the base case. His proof is inspired by the following property of diffeomorphism groups which is called fragmentation.
\begin{prop*} Let $M$ be a compact manifold and  let $\{ U_i\}_{i}$ be a finite open cover of $M$. Then any element $f\in \Diff^r_0(M)$ can be written as a composition of diffeomorphisms $f_j$ such that  $f_j$ is compactly supported in some element of the cover  $\{ U_i\}_{i}$.
\end{prop*}
To illustrate the idea, let us focus on section spaces of fiber bundles and how Thurston fragmented them. Let $\pi\colon E\to M$ be a Serre fibration over an $n$-dimensional Riemannian manifold $M$ with a non-zero injectivity radius. Let  $s_0$ be the fixed section of this bundle as a base section. So with respect to this base section, we define $\text{Sect}_c(\pi)$ to be the space of compactly supported sections of $\pi$ with compact open topology. 

We want to give a filtration on $\text{Sect}_c(\pi)$ that has a ``nice'' filtration quotients. Fix a positive $\epsilon$ smaller than the injectivity radius. Let $\text{Sect}_{\epsilon}(\pi)$ denote the subspace of sections $s$ such that the support of $s$ can be covered by $k$ geodesically convex balls of radius $2^{-k}\epsilon$ for some positive integer $k$; we call this subset $\epsilon$-supported sections. 

There is a filtration on  $\text{Sect}_{\epsilon}(\pi)$ whose $k$-th level is the subset of sections whose supports can be covered by at most $k$ balls. And the reason for the choice of $2^{-k}\epsilon$ is that the filtration quotients behave nicely. For example, if two balls of radius $\epsilon/4$ intersect, then they can be contained in a ball of radius $\epsilon/2$. So a section in the second filtration quotient has support in two disjoint balls. This is homotopically useful for studying filtration quotients.

 Thurston proved the following which can be improved (see \cite{MR4580303}) to prove what is now called {\it non-abelian Poincar\' e duality}.
\begin{thm}[Thurston's fragmentation] If the fiber of $\pi$ is at least $n$-connected, the inclusion 
\[
\text{\textnormal{Sect}}_{\epsilon}(\pi)\hookrightarrow \text{\textnormal{Sect}}_c(\pi)
\]
is a weak homotopy equivalence. 
\end{thm}
Since the fiber of $\pi$ is $n$-connected, the space $\text{\textnormal{Sect}}_c(\pi)$ is connected so each section can be deformed to the base section, which lies in $\text{\textnormal{Sect}}_{\epsilon}(\pi)$. To prove the theorem, one has to deform a family of compactly supported sections to a family of $\epsilon$-supported sections. To see the main idea, let us see how we can obtain such a deformation for a one-parameter family of sections in $\text{\textnormal{Sect}}_c(\pi)$.

 Let $f\colon [0,1]\to \text{\textnormal{Sect}}_c(\pi)$ be a one-parameter family. Consider the adjoint of this map $F\colon [0,1]\times M\to E$. We homotope $F$ in two steps so that at each time coordinate we get an $\epsilon$-supported section.  

Let $\{\mu_i\}_{i=1}^N$ be a partition of unity with respect to an open cover of $M$. We define a {\it fragmentation homotopy} with respect to this partition of unity. Let $\nu_j$ be the function
\[
\nu_j(x)=\sum_{i=1}^j\mu_i(x).
\]

\[
H_1:[0,1]\times M\to [0,1]\times M,
\]
\[
H_1(t,x)=(u,x),
\]
\[
u(t,x)=\nu_{\floor{Nt}}(x)+\mu_{\floor{Nt}+1}(x)(Nt-\floor{Nt}).
\]
Since $H_1(x,t)$ preserves the $x$ coordinate, we can define a straight line homotopy $H_t\colon [0,1]\times M\to [0,1]\times M$ from the identity to $H_1$. 
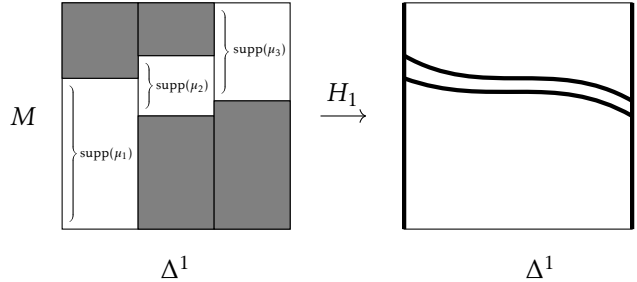
\begin{figure}[h]\label{frag}
\[
\begin{tikzpicture}
  \def\rectanglepath{-- ++(3cm,0cm)  -- ++(0cm,3cm)  -- ++(-3cm,0cm) -- cycle}
  \draw (0,0) \rectanglepath;
   \draw (4.5,0) \rectanglepath;

  \draw [decorate,
    decoration = {calligraphic brace, mirror}] (0.1,0.1) --  (0.1,1.9);
 
 \draw [decorate,
    decoration = {calligraphic brace, mirror}] (1.1,1.6) --  (1.1,2.2);
    
     \draw [decorate,
    decoration = {calligraphic brace, mirror}] (2.1,1.8) --  (2.1,2.9);
 
\draw (1,0)--(1,3);
\draw (2,0)--(2,3);
\draw (0,2)--(1,2);
\draw (1, 1.5)--(2, 1.5);
\draw (1,2.3)--(2, 2.3);
\draw (2, 1.7)--(3,1.7);
\node (X) at (1.5,-0.5) {$\Delta^1$};
\node (Y) at (6.3,-0.5) {$\Delta^1$};
\node (Z) at (-0.5,1.5) {$M$};
\node (R) at (0.58,1.) {\scalebox{0.5}{$\text{supp}(\mu_1)$}};
\node (R) at (1.58,1.9) {\scalebox{0.5}{$\text{supp}(\mu_2)$}};
\node (R) at (2.6,2.37) {\scalebox{0.5}{$\text{supp}(\mu_3)$}};
\draw [ fill=gray] (0,2) rectangle (1,3);
\draw [fill=gray] (1,0) rectangle (2,1.5);
\draw [ fill=gray] (1,2.3) rectangle (2,3);
\draw [ fill=gray] (2,0) rectangle (3,1.7);
\draw [ultra thick] (4.5,2) to [out=-20,in=150] (7.5,1.5);
\draw [ultra thick] (4.5,2.3) to [out=-30,in=150] (7.5,1.7);
\draw [ultra thick] (4.5,0) -- (4.5,3);
\draw [ultra thick] (7.5,0) -- (7.5,3);
\draw [->](3.4,1.5) -- (4.,1.5) node[midway,above] {$H_1$};
\end{tikzpicture}
\]

\caption{Fragmentation map for $N=3$ and $q=1$. The thick lines are the images of $\{0\}\times M, \{1/3\}\times M, \{2/3\}\times M$ and $ \{1\}\times M$ under the map $H_1$.}
\end{figure}
As in Figure \ref{frag}, the map $H_1$ is defined so that the gray area is mapped to a subcomplex $L$ of  $[0,1]\times M$ (the bold lines in Figure \ref{frag}) in the target which is of dimension $n=\text{dim}(M)$. 

Using the fact that $\overline{\mathrm{B}\Gamma}_n^r$  is at least $n$-connected and $L$ is $n$-dimensional subcomplex, it is standard by the obstruction theory, to change $F\colon [0,1]\times M\to E$ up to homotopy to a map $G\colon [0,1]\times M\to E$ such that the restriction of $G$ to the subcomplex $L$ is the same as $\text{id}\times s_0$. Then $F$ is homotopic to $ G\circ H_1$ and one can see that $ G\circ H_1\colon [0,1]\times M\to E$ gives a one-parameter family that at each time is $\epsilon$-supported section. This idea also works for higher dimensional parameter space to prove fragmentation for section spaces. Note that in the Mather-Thurston theorem the fiber of the bundle $\tau_M^*(\nu^r)$ is $\overline{\mathrm{B}\Gamma}_n^r$ which is at least $n$-connected so the above theorem applies and we obtain a filtration on the section space.  

On the other hand, recall that the support for foliated trivialized bundle $B\times M\to B$ is the maximal closed subset $K$ of the fiber $M$ such the foliation on $B\times (M\backslash K)\to B$ is horizontal foliation. And one can similarly prove a fragmentation property for the classifying spaces that classify foliated trivialized bundles to obtain a filtration. One can also fragment the Haefliger space $\overline{\mathrm{B}\Gamma}_n^r$ which classifies germs foliations by ``thickening'' its model as follows. The set of $k$-simplices of the semi-simplicial set that realizes to $\overline{\mathrm{B}\Gamma}_n^r$ is the set of germs of foliated bundle $\Delta^k\times \bR^n\to \Delta^k$ near the zero section $\Delta^k\times \{0\}$. We can thicken this model to define the set of $k$-simplices to be the set of germs of foliated bundle $\Delta^k\times \bR^n\to \Delta^k$ near the disk bundle $\Delta^k\times D^n$. This model is also homotopy equivalent to also $\overline{\mathrm{B}\Gamma}_n^r$ and now one considers the support as a subset of the disk $D^n$ and the fragmentation idea applies. Thurston considered the  map
\begin{equation}\label{N}
f\colon\overline{\mathrm{B}\Gamma}_n^r\xrightarrow{\simeq} \text{Map}(D^n, \overline{\mathrm{B}\Gamma}_n^r),
\end{equation}
that is induced by constant maps. So now both $\overline{\mathrm{B}\Gamma}_n^r$ and $\text{Map}(D^n, \overline{\mathrm{B}\Gamma}_n^r)$ satisfy fragmentation property and the map between them is filtration preserving. 
\[
\begin{tikzcd}
X_1 \arrow{d}{f_1}\arrow[r,hook] &X_2 \arrow{d}{f_2}\arrow[r,hook] &\cdots \arrow[r,hook] & \overline{\mathrm{B}\Gamma}_n^r \arrow{d}{\simeq}\\
Y_1\arrow[r,hook] &Y_2 \arrow[r,hook] &\cdots \arrow[r,hook] & \text{Map}(D^n, \overline{\mathrm{B}\Gamma}_n^r).
\end{tikzcd}
\]
It turns out that $X_1$ and $Y_1$ are homology isomorphic to the $n$-fold suspensions $\Sigma^n(\overline{\BDiff_c^r(\bR^n)}) $ and $\Sigma^n(\Omega^n \overline{\mathrm{B}\Gamma}_n^r)$ respectively and proving that $f_1$ induces a homology isomorphism implies that 
\[
\Sigma^n(\overline{\BDiff_c^r(\bR^n)})\to \Sigma^n(\Omega^n \overline{\mathrm{B}\Gamma}_n^r),
\]
also induces a homology isomorphism which implies the same for the map \ref{T}. Because suspension only shifts the homological degree.

 The key is to observe that the map between filtration quotients
\[
\overline{f_k}\colon X_k/X_{k-1}\to Y_k/Y_{k-1},
\]
induces a homology isomorphism in a range that increases linearly in $k$ and the fact that in the limit we have the weak equivalence \ref{N}, a standard argument of comparison of spectral sequences implies that $f_1$ induces a homology isomorphism.

\section{The Equivariant Mather-Thurston theorem} Recall that for a topological group $G$, the universal principal $G$-bundle is $\mathrm{E}G\to\mathrm{B}G$, where $\mathrm{E}G$ is a contractible space on which $G$ acts freely. Given a map $f\colon X\to \mathrm{B}G$, one can choose a model for $F$, the homotopy fiber of the map $f$, on which $G$ acts. For example, the pullback $f^*(\mathrm{E}G)$ is a model for the homotopy fiber of $f$ on which $G$ acts naturally. It is standard to see that the homotopy quotient of this action
\[
f^*(\mathrm{E}G)\hcoker G\coloneqq \frac{f^*(\mathrm{E}G)\times \mathrm{E}G}{G},
\]
is homotopy equivalent to $X$. Therefore, we can choose a model for $\overline{\BDiff_c^r(M)}$, which is the homotopy fiber of 
\[
\eta\colon \BcdrDiff M\to \BDiff^r_c(M),
\]
on which $\Diff_c^r(M)$ acts and the homotopy quotient $\overline{\BDiff_c^r(M)}\hcoker \Diff_c^r(M)$ is weakly equivalent to $ \BcdrDiff M$.

There is also a natural model for $\text{Sect}_c(\tau_M^*(\nu^r))$ in the Mather-Thurston theorem \ref{MT} on which $\Diff_c^r(M)$ acts as follows. Let $\gamma_n$ be the tautological $n$-dimensional vector bundle over $\mathrm{B}\text{GL}_n(\bR)$. Let $\mathrm{Bun}(\mathrm{T}M, (\nu^r)*\gamma_n)$  be the space of bundle maps $\mathrm{T}M\to (\nu^r)^*\gamma_n$ i.e. continuous maps that are linear isomorphisms on each fiber, equipped with the compact open topology. It is easy to define a map
\[
\text{Sect}_c(\tau_M^*(\nu^r))\to \mathrm{Bun}(\mathrm{T}M, (\nu^r)^*\gamma_n),
\]
which is a weak homotopy equivalence. The group $\Diff_c^r(M)$ acts on $\mathrm{Bun}(\mathrm{T}M, (\nu^r)^*\gamma_n)$ by precomposition with the differential of diffeomorphisms.  In my thesis, I observed that 
 \begin{thm}
  The map in the Mather-Thurston theorem \ref{MT} can be made  $\Diff_c^r(M)$-equivariant.
 \end{thm}  
 Hence, the classifying space $ \BcdrDiff M$ is homology isomorphic to the space $\mathrm{Bun}(\mathrm{T}M, (\nu^r)^*\gamma_n)\hcoker \Diff_c^r(M)$ which I used  to study the stable homology of surface diffeomorphism groups. Let me mention one consequence of this equivariance about invariants of flat bundles.
 
 Recall that oriented circle bundles are classified by their Euler class. And we have $\BDiff_0^r(S^1)\simeq \mathrm{B}S^1$, whose integral cohomology is a polynomial ring on the Euler class. Hence, we have a map
 \[
 \bQ[e]\cong H^*(\mathrm{B}S^1; \bQ)\to H^*(\BdrDiffo {S^1}; \bQ).
 \]
 \begin{thm}[Morita]The powers of the Euler class $e^k$ are nontrivial in $H^*(\BdrDiffo {S^1}; \bQ)$ for all $k$.
 \end{thm}
 If we restrict the holonomy of the bundle to $\text{PSL}_2(\bR)$, this is no longer true and it turns out that $e^2$ vanishes in $H^4(\mathrm{B}\text{PSL}_2(\bR)^{\delta};\bQ)$. Morita used rational homotopy theory to prove this theorem but one can give a simpler proof using the equivariant the Mather-Thurston theorem. The space $\BdrDiffo {S^1}$ is homology isomorphic to $\mathrm{Bun}(\mathrm{T}S^1, (\nu^r)^*\gamma_1)\hcoker \Diff_0^r(S^1)$. Since $\Diff_0^r(S^1)$ acts nontrivially on the tangent vectors of $S^1$, the action on $\mathrm{Bun}(\mathrm{T}S^1, (\nu^r)^*\gamma_1)$ {\it does not} have a fixed point. But  $\Diff_0^r(S^1)$ is homotopy equivalent to the subgroup of rotations $S^1\subset \Diff_0^r(S^1)$. And it turns out that the action of $S^1$ on $\mathrm{Bun}(\mathrm{T}S^1, (\nu^r)^*\gamma_1)$ has fixed points. This implies that 
 \[
 \mathrm{Bun}(\mathrm{T}S^1, (\nu^r)^*\gamma_1)\hcoker S^1\to \mathrm{B}S^1,
 \]
 has a section that gives a simpler proof of Morita's theorem. 
 \begin{problem}Let $\Diff^{\omega}_0(S^1)$ denote the orientation preserving analytic diffeomorphisms of $S^1$. Are the powers of the Euler class $e^k$ nontrivial in $H^*(\BDiff_0^{\omega}(S^1)^{\delta};\bQ)$?
 \end{problem}
 Thurston claimed that $e^3$ vanishes for analytic flat circle bundles but the proof was not written down and Etienne Ghys called this claim the lost theorem of Thurston!
 
 In the above fixed-point argument, it is important that $S^1$ is a Lie group. For higher dimensional sphere bundles one can use torus actions to prove a generalization of Morita's theorem.

\begin{thm}[Nariman]
The classes $e^k$ for all $k$  are all nontrivial in $H^*(\BdrDiffo {S^{2n-1}};\bQ)$.
\end{thm}

\section{Controlled Mather-Thurston theorems and Milnor-Wood inequalities}\label{CMT} Recall that the Mather-Thurston theorem for $C^0$-foliations and Mather's acyclicity result (\cite{MR0288777}) imply that 
\[
\BdH_0(M)\to \BH_0(M),
\]
induces a homology isomorphism. Geometrically, this means for each $M$-bundle $M\to E\xrightarrow{\pi} B$ whose group structure is $\tH_0(M)$, there exists a  bordism $V$ from $B$ to $B'$ and an $M$ bundle $W\to V$ that extends the bundle $\pi$ and its restriction to $B'$ is a flat $M$-bundle. Intuitively, one expects that the bordism $V$ from $B$ to $B'$ ought to make the fundamental group of $B$ more and more complicated until $B'$ could accommodate a flat bundle with the same characteristic numbers as the bundle $\pi$. 

For example for $M=S^1$, a topological oriented $S^1$-bundle $E\xrightarrow{p} \Sigma_g$ over $\Sigma_g$ a closed oriented surface of genus $g$ is classified by its Euler class $e(p)$. The Milnor-Wood inequality in this case gives the necessary and sufficient condition under which the circle bundle $p$ admits a flat structure, namely
\[
|\langle e(p), [\Sigma_g]\rangle|\leq 2g-2,
\]
where the left-hand side is the evaluation of the class $e(p)$ on the fundamental cycle $[\Sigma_g]$ and it is called the Euler number of the circle bundle $p$. In a more modern language, the Euler class for flat $S^1$-bundles is bounded in the sense of Gromov. So fixing the Euler number of a circle bundle, in order to make the bundle flat up to bordism, we have to increase the genus of the base so that it satisfies the Milnor-Wood inequality. 

Inspired by physics, Mike Freedman in \cite{freedman2020controlled} took a constructive point of view and demanded two types of control on the bordism $V$ whose existence is guaranteed by the Mather-Thurston theorem. One of them is to formalize the idea that the bordism $V$ is directed from $B$ to $B'$ to increase ``the complexity'' of $\pi_1(B)$ and the other is to control the holonomy of the bundle on $B'$. 

\begin{defn}
A  semi-s-cobordism  is a manifold triple $(V, B, B')$ with $\partial V=B\coprod -B'$ where $-B'$ means $B'$ with the reverse orientation, so that the inclusion $\iota\colon B\hookrightarrow V$ is a ``{\it simple}'' homotopy equivalence.
\end{defn} 
The directionality of $V$ in the definition is that we do not assume that $\iota'\colon B'\hookrightarrow V$ is a homotopy equivalence. However, the definition implies that there is a (simple) deformation retraction $r\colon V\to B$ and the map $r\circ \iota'\colon B'\to B$ is a degree one map which implies that the map $\pi_1(B')\to \pi_1(B)$ is surjective. 

Freedman uses low dimensional techniques to construct certain ``homological'' solid tori and the Bing double construction to do surgeries to construct the bordism $V$ in dimension $4$.
\begin{thm}[Freedman \cite{freedman2020controlled}] Suppose $B$ is a closed $3$-manifold and $M\to E\xrightarrow{\pi} B$ is an $M$-bundle. Then there exists a  {\it semi-s-cobordism} $V$ from $B$ to $B'$ and an $M$ bundle $W\to V$ that extends the bundle $p$ and its restriction to $B'$ is a flat $M$-bundle.
\end{thm}
 Freedman conjectured that not only this statement holds in all dimensions but also even more strongly when $M$ is a Riemannian manifold, for every positive $\epsilon$, one can choose the semi-s-cobordism $V$ so that the holonomy of the flat $M$ bundle on $B'$ on a generating set of $\pi_1(B')$ is at most $\epsilon$ away from the isometry group of $M$ in the sup norm on the group $\tH_0(M)$. 

 In fact in higher dimensions, we later learned that the existence of such semi-s-cobordism was known using homological surgery techniques.
 \begin{thm}[Hausmann and Vogel \cite{MR0511783}] Let $M$ and $B$ be two closed topological manifolds where $\text{dim}(B)\geq 5$ and let $\xi$ be an $M$-bundle over $B$. Then this $\xi$ extends as an $ M$ bundle over $B$ to a semi-s-cobordism such that at the other end, it is a $C^0$-foliated $ M$ bundle. 
 \end{thm}

Demanding the bordism to be semi-s-cobordism is a qualitative control and the control on the holonomy which is related to the Milnor-Wood inequality is more quantitative. Given the existence of semi-s-cobordism reduction to $C^0$-foliated bundles is known in all dimensions except $4$, one interesting unknown case is the following.

\begin{problem}
Suppose $S^3\to S^7\to S^4$ is a ``generalized'' Hopf fibration. Does there exist a homology $4$-sphere $H$ and a degree one map $H\to S^4$ such that the pull-back of the generalized Hopf fibration is a flat $S^3$-bundle?
\end{problem}
In fact in the 80s, Ghys asked if the Milnor-Wood inequality holds for oriented flat $S^3$-bundles.
\begin{quest}[Ghys]
Let $M^4$ be a compact orientable $4$-manifold and $ \pi_1(M) \xrightarrow{\rho} \Homeo_\circ(\bS^3)$ be a representation. Is it true that the Euler number of the flat $S^3$-bundle associated with $\rho$ over $M$ is bounded by a number depending only on $M$?
\end{quest}
 Monod and I answered Ghys' question negatively (\cite{MR4588567}) by proving the following result.
 \begin{thm}(\cite[Theorem 1.8]{MR4588567})
The fourth bounded cohomology of the group $\tdH_0(S^3)$  vanishes.
\end{thm}
This, in particular, implies that all non-trivial classes in $H^4(\BdH_0(S^3);\bR)$ are unbounded so the Euler class for oriented flat $S^3$ bundle is not bounded. Hence, in Freedman's conjecture, we might have the semi-s-cobordism to flat bundles in higher dimensions, but we may not be able to quantitatively control the topology of $B'$. However, the unboundedness of the Euler class seems to suggest that one might be able to flatten $S^3$-bundles over $4$-manifolds up to bordism without changing the simplicial volume of the base much.

Our proof of the failure of the Milnor-Wood inequality for $S^3$-bundles is not constructive and it leads to both conceptual and computational questions (\cite[Section 7]{MR4588567}) but most importantly we lack geometric intuition to construct examples to exhibit the failure of Milnor-Wood in this case.

 \subsection*{Acknowledgments} I am grateful to Laszlo Lempert, Mehdi Yazdi, Fran\c cois Laudenbach, and Mehrdad Shahshahani for their helpful comments on the first draft of this article to improve its readability.  The author was partially supported by  NSF CAREER Grant DMS-2239106 and Simons Foundation Collaboration Grant (855209). I would like to thank Allison May for the opening figure, Patrick Massot for the figure $8$, and my father for helping me with other figures.

\bibliography{bib}
\bibliographystyle{alpha}

\end{document}